# New Results on Sampled-Data Feedback Stabilization for Autonomous Nonlinear Systems

J. Tsinias
Department of Mathematics, National Technical University
Zografou Campus 15780, Athens, Greece
jtsin@central.ntua.gr

**Abstract.** Sufficient conditions are established for sampled-data feedback global asymptotic stabilization for nonlinear autonomous systems. One of our main results is an extension of the well known Artstein-Sontag theorem on feedback stabilization concerning affine in the control systems. A second aim of the present work is to provide sufficient conditions for sampled-data feedback asymptotic stabilization for two interconnected nonlinear systems. Lie algebraic sufficient conditions are derived for the case of affine in the control interconnected systems without drift terms.

## 1. Introduction

In the recent literature on nonlinear control theory an important area that has received much attention is the stabilization problem by means of sampled-data feedback; see for instance, [1],[6]-[8],[10],[11],[16]-[27] and relative references therein, where sufficient conditions are established for the existence of sampled-data and hybrid feedback controllers exhibiting stabilization. We also mention the recent contributions [3],[4],[21] and [30], where, under the presence of Control Lyapunov Functions (CLF), "triggering" techniques are developed for the determination of the set of the sampling time instants for the corresponding sampled-data controller. Particularly, in [4] a "minimum attention control" approach is adopted, exhibiting minimization of the open loop operation of the sampled-data control. In [21] a universal formula is proposed for the sample-data feedback exhibiting "event-based" stabilization of affine in the control systems. The corresponding result constitutes a generalization of the well known Sontag's result in [29].

In the recent author's works [32]-[34] the concept of *Weak Global Asymptotic Stabilization by Sampled-Data Feedback* (SDF-WGAS) is introduced and Lyapunov-like sufficient conditions for the existence of a sampled-data feedback stabilizer have been established. These conditions are weaker than those proposed in earlier contributions on the same subject. The present paper constitutes a continuation of the previously mentioned author's works. The paper is organized as follows:

The current section contains the precise definition of the SDF-WGAS, as originally given in [33]. In Section 2 we establish a general result (Proposition 1), which provides Lyapunov characterizations of certain properties relying asymptotic controllability for general nonlinear systems:

$$\dot{x} = F(x,u), (x,u) \in \mathbb{R}^n \times \mathbb{R}^\ell \tag{1.1a}$$

$$F(0,0) = 0 \tag{1.1b}$$

where we assume that $F : \mathbb{R}^n \times \mathbb{R}^m \to \mathbb{R}^n$ is Lipschitz continuous.

Sections 3 and 4 are devoted to applications of Proposition 1 for the solvability of the SDF-WGAS problem for certain class of nonlinear systems. In Section 3 we apply the result of Proposition 1 to derive an extension of the so called Artstein-Sontag theorem (see [5], [29] and also [6],[7] for recent extensions). Particularly, in Proposition 2 a sufficient condition for the solvability of the SDF-WGAS problem for the case of affine in the control systems:

$$\dot{x} = f(x) + ug(x), (x,u) \in \mathbb{R}^n \times \mathbb{R} \tag{1.2a}$$

$$f(0) = 0 \tag{1.2b}$$

is established. This condition is weaker than the familiar hypothesis imposed in [5] and includes the Lie bracket between the vector fields $f$ and $g$. In Section 4 we use the result of Proposition 1 to provide a small-gain criterion for the possibility of sampled-data feedback global stabilization for composite systems of the form:

$$\Sigma_1 : \dot{x} = f(x, y, u)$$
$$\Sigma_2 : \dot{y} = g(x, y, u) \tag{1.3a}$$
$$(x, y, u) \in \mathbb{R}^n \times \mathbb{R}^m \times \mathbb{R}^\ell$$
$$f(0,0,0) = 0, g(0,0,0) = 0 \tag{1.3b}$$

The corresponding result (Proposition 3) extends the main result in [34] and partially extends the well known results in the literature (see for instance [14],[15]) establishing small-gain criteria for composite systems $\Sigma_1$, $\Sigma_2$ with no controls and particularly those in [2],[12],[13] (see also relative references therein), where each subsystem does not necessarily satisfies the Input-to-State-Stability (ISS) property. As a consequence of Proposition 3, we provide a partial extension of the well known result in [9] due to J.-M.Coron, concerning the solvability of the stabilization problem by means of smooth time-varying feedback for the affine in the control systems without drift term. In the present work we consider nonholonomic composite systems of the following form:

$$\dot{\xi} = F(\xi, u) := \sum_{i=1}^{\ell} u_i F_i(\xi) \tag{1.4a}$$

$$F_i(\xi) := \begin{pmatrix} A_i(\xi) \\ B_i(\xi) \end{pmatrix}, i = 1, ..., \ell, \ \xi := (x, y) \in \mathbb{R}^n \times \mathbb{R}^m \tag{1.4b}$$

$$A_i(0) = 0, B_i(0) = 0, i = 1, ..., \ell \tag{1.4c}$$

and in Proposition 4 we derive a set of Lie algebraic sufficient conditions guaranteeing SDF-WGAS for the case (1.4), being weaker than the *accessibility rank condition* imposed in [9].

*Notation and Definitions:* Throughout the paper, we adopt the following notation. By $x^T$ we denote the transpose of a given vector $x \in \mathbb{R}^n$. By $K$ we denote the set containing all continuous strictly increasing functions $\phi: \mathbb{R}^+ \to \mathbb{R}^+$ with $\phi(0) = 0$ and $K_\infty$ denotes the subset of $K$ that is constituted by all $\phi \in K$ with $\phi(t) \to \infty$ as $t \to \infty$. We denote by $\pi(\cdot, s, x_0, u)$ the trajectory of (1.1a) with $\pi(s, s, x_0, u) = x_0$ corresponding to certain (measurable and essentially bounded) control $u : [s, T_{\max}) \to \mathbb{R}^m$, where $T_{\max}$ is the corresponding maximal existing time of the trajectory.

**Definition 1.** We say that (1.1) is *Globally Asymptotically Controllable* at zero (GAC), if for any $x_0 \in \mathbb{R}^n$ there exists a control $u_{x_0}(\cdot) : \mathbb{R}^+ \to \mathbb{R}^\ell$ such that $\pi(t, 0, x_0; u_{x_0})$ exists for all $t \geq 0$ and the following properties hold:

*Stability:* $\quad \forall \varepsilon > 0 \Rightarrow \exists \delta = \delta(\varepsilon) > 0 : |x_0| \leq \delta \Rightarrow |\pi(t, 0, x_0, u_{x_0})| \leq \varepsilon, \ \forall t \geq 0 \tag{1.5a}$

*Attractivity:* $\quad \lim_{t \to \infty} \pi(t, 0, x_0, u_{x_0}) = 0, \ \forall x_0 \in \mathbb{R}^n \tag{1.5b}$

**Definition 2.** We say that (1.1) is *Weakly Globally Asymptotically Stabilizable by Sampled-Data Feedback* (SDF-WGAS), if for any constant $\sigma > 0$ there exist a map $T : \mathbb{R}^n \to \mathbb{R}^+$ satisfying $T(0) = 0$ and

$$T(x) \leq \sigma, \quad \forall x \in \mathbb{R}^n \setminus \{0\} \tag{1.6}$$

and a map $\varphi := \varphi(t, s, x) : \mathbb{R}^+ \times \mathbb{R}^+ \times \mathbb{R}^n \to \mathbb{R}^\ell$ with $\varphi(\cdot, \cdot, 0) = 0$, being measurable and essentially bounded with respect to the first two variables, such that for every $x_0 \neq 0$, a sequence of times

$$t_0 := 0 < t_1 < t_2 < ... < t_\nu < ...; \quad t_\nu \to \infty \tag{1.7}$$

can be found, in such a way that, if we denote:

$$u_{x_0}(t) := \varphi(t, t_i, \pi(t_i)), \quad t \in [t_i, t_{i+1}) \quad i = 0, 1, 2, \ldots \tag{1.8a}$$

$$\pi(t_i) := \pi(t_i, t_{i-1}, \pi(t_{i-1}), u_{x_0}), \quad i = 1, 2, \ldots; \pi(t_0) := x_0 \tag{1.8b}$$

then

$$t_{n+1} - t_n = T(\pi(t_n)), \quad n = 0, 1, 2, \ldots \tag{1.9}$$

and (1.1) is GAC by means of the controller $u_{x_0}(\cdot)$ as defined in (1.8); equivalently, the system

$$\dot{x}(t) = f(x(t), u_{x_0}(t)), \quad t \in [t_i, t_{i+1})$$
$$\text{with } u_{x_0}(\cdot) \text{ as defined in (1.8)} \tag{1.10}$$

satisfies properties (1.5a) and (1.5b).

## 2. A General Result

Consider the system (1.1) and assume that there exist a closed set $D_1 \subset \mathbb{R}^n$ containing zero $0 \in \mathbb{R}^n$, an open neighborhood $D_2$ of $D_1$, constants

$$0 \le L_1 < L_2 \le +\infty, \tag{2.1}$$

a continuous function $\Phi : \mathbb{R}^n \to \mathbb{R}^+$ with

$$\sup\{\Phi(x) : x \in D_1\} \le L_1; \tag{2.2a}$$

$$0 < \Phi(x) < L_2, \forall x \in D_2 \setminus D_1; \tag{2.2b}$$

$$\Phi(x) \ge L_2, \forall x \in \mathbb{R}^n \setminus D_2, \text{ provided that } \mathbb{R}^n \setminus D_2 \ne \emptyset; \tag{2.2c}$$

$$\Phi(x) = L_2, \forall x \in \partial D_2 \text{ and for every sequence } \{x_\nu \in D_2 \setminus D_1\} \text{ with } \lim_{\nu \to \infty} |x_\nu| = \infty \text{ it holds}$$

$$\limsup_{\nu \to \infty} \Phi(x_\nu) = L_2 \tag{2.2e}$$

and a function $\beta \in K$, such that for every constant $\sigma > 0$ and for every $x \in D_2 \setminus D_1$ there exists a time $\tau = \tau(x) \in (0, \sigma]$ and a control $u_x(\cdot) : [0, \tau] \to \mathbb{R}^m$ satisfying the following properties:

$$\Phi(\pi(\tau, 0, x, u_x)) < \Phi(x), \tag{2.3a}$$

$$\Phi(\pi(t, 0, x, u_x)) \le \beta(\Phi(x)), \forall t \in [0, \tau] \tag{2.3b}$$

The first two statements of the following proposition generalize the result in [33, Proposition 1] and are used in Section 4 for the derivation of sufficient conditions for SDF-WGAS for the cases (1.3) and (1.4).

**Proposition 1.** For system (1.1) assume that (2.1)-(2.3) are fulfilled. Then the following hold:
(i) For any constant $\sigma > 0$ there exists a map $T : D_2 \setminus D_1 \to \mathbb{R}^+ \setminus \{0\}$ satisfying (1.6) and a map $\varphi := \varphi(t, s, x) : \mathbb{R}^+ \times \mathbb{R}^+ \times (D_2 \setminus D_1) \to \mathbb{R}^\ell$ with same regularities properties with those in Definition 2, such that for every $x_0 \in D_2 \setminus D_1$ there exists an increasing sequence of times $\{t_\nu\}$, in such a way that, if we denote $u_{x_0}(t) := \varphi(t, t_i, \pi(t_i))$, $t \in [t_i, t_{i+1}), i = 0, 1, 2, \ldots$, $\pi(t_i) := \pi(t_i, t_{i-1}, \pi(t_{i-1}), u_{x_0})$, $i = 1, 2, \ldots$; $\pi(t_0) := x_0$, then (1.9) holds and, if we consider the resulting system (1.10), then for every neighborhood $N \subset D_2$ of $D_1$ there exists a time $\tau := \tau(x_0) > 0$ such that its trajectory initiated from $x_0 \in D_2 \setminus N$ satisfies

$$\pi(\tau, 0, x_0, u_{x_0}) \in \text{int } N \tag{2.4}$$

(ii) If, in addition to previous assumptions, we assume that

$$D_1 = \{0\}, L_1 = 0 \tag{2.5}$$

then for the trajectories of (1.10) both properties (1.5a,b) hold, provided that $x_0 \in D_2$.
(iii) If, in addition to (2.5), we assume that

$$D_2 = \mathbb{R}^n, L_2 = \infty \tag{2.6}$$

then system (1.1) is SDF-WGAS.

**Remark 1.** (i) The last statement of Proposition 1 coincides with the main result obtained in [33] and constitutes a generalization of Theorem 17 in [28]. To be precise, Proposition 1 in [33] asserts that system (1.1) is SDF-WGAS, if there exist a continuous, positive definite and proper function $\Phi : \mathbb{R}^n \to \mathbb{R}^+$ and a function $\beta \in K$ such that, for every constant $\sigma > 0$ and for every $x \in \mathbb{R}^n \setminus \{0\}$, a time $\tau = \tau(x) \in (0, \sigma]$ and a control $u_x(\cdot) : [0, \tau] \to \mathbb{R}^\ell$ can be found satisfying:

$$\Phi(\pi(\tau, 0, x, u_x)) < \Phi(x), \tag{2.7a}$$

$$\Phi(\pi(t, 0, x, u_x)) \leq \beta(\Phi(x)), \forall t \in [0, \tau] \tag{2.7b}$$

Notice that both (2.7a,b) are fulfilled with $\beta(s) = s$, if we assume that for every constant $\sigma > 0$ and for every $x \in \mathbb{R}^n \setminus \{0\}$ there exists a time $\tau = \tau(x) \in (0, \sigma]$ and a control $u_x(\cdot) : [0, \tau] \to \mathbb{R}^\ell$ satisfying:

$$\Phi(\pi(\tau, 0, x, u_x)) < \Phi(x), \forall t \in (0, \tau] \tag{2.8}$$

Obviously, (2.8) is fulfilled, if $\Phi$ is a $C^1$ (global) CLF for (1.1), namely, if we assume that for every $x \in \mathbb{R}^n \setminus \{0\}$ there exists a vector $u \in \mathbb{R}^m$ satisfying $D\Phi(x)F(x,u) < 0$. Of course, the converse claim is not in general true. It turns out, that the above hypotheses (2.7a,b) guaranteeing SDF-WGAS, are weaker than those adopted in earlier works on the literature, that mainly based on the existence of a $C^1$ CLF.

(ii) It is worthwhile noticing that, as is shown in [33, Proposition 1], the statement (iii) of Proposition 1 is equivalent to the existence of mappings $\Phi : \mathbb{R}^n \to \mathbb{R}^+, a_1, a_2 \in K_\infty$, $\Phi$ being generally *discontinuous* on $\mathbb{R}^n \setminus \{0\}$ with $a_1(|x|) \leq \Phi(x) \leq a_2(|x|), \forall x \in \mathbb{R}^n$ and a function $\beta \in K$ such that for every constant $\sigma > 0$, a pair of lower semi-continuous (LSC) mappings $T : \mathbb{R}^n \setminus \{0\} \to \mathbb{R}^+ \setminus \{0\}$ and $q : \mathbb{R}^n \to \mathbb{R}^+ \setminus \{0\}$ can be found in such a way that (1.6) holds and for every $x \in \mathbb{R}^n \setminus \{0\}$ there exists a control $u_x(\cdot) : [0, T(x)] \to \mathbb{R}^\ell$ satisfying $\Phi(\pi(T(x), 0, x, u_x)) \leq \Phi(x) - q(x)$ and $\Phi(\pi(t, 0, x, u_x)) \leq \beta(\Phi(x)), \forall t \in [0, T(x)]$.

**Proof of Proposition 1:** From (2.3) and by taking into account continuity of the map $\Phi$ and continuity of the solution of (1.1) with respect to the initial values, it follows that for every $x \in D_2 \setminus D_1$, a time $\tau = \tau(x) \in (0, \sigma]$ and a control $u = u_x : [0, \tau] \to \mathbb{R}^\ell$ can be found such that, if we denote $\bar{\beta} := 2\beta$, the following holds:

$$\Phi(\pi(\tau, 0, y, u)) < \Phi(y) - q, \tag{2.9a}$$

$$\Phi(\pi(t, 0, y, u)) \leq \bar{\beta}(\Phi(y)), \forall t \in [0, \tau], \text{ for every y near } x \tag{2.9b}$$

for certain constant $q := q(x) > 0$. Using (2.9) we can determine a locally finite covering of $D_2 \setminus D_1$ of closed spheres $S_\nu := S[x_\nu, r_\nu] \subset D_2 \setminus D_1$, $\nu = 1, 2...$ of radius $r_\nu$ centered at $x_\nu \in D_2 \setminus D_1$, times $\tau_\nu$, controls $u_\nu : [0, \tau_\nu] \to \mathbb{R}^\ell$ and constants $q_\nu > 0$, in such a way that both (2.9a,b) are fulfilled, with $x := x_\nu, \tau := \tau_\nu$, $u := u_\nu$ and $q := q_\nu$, for every $y \in S_\nu, t \in [0, \tau_\nu]$. For each $x \in D_2 \setminus D_1$ let $\{S_{\nu_i}, i = 1, ...., j\}$ be the sequence of those spheres above for which $x \in \bigcap S_{\nu_i}$ for $1 \leq i \leq j$ and $x \notin S_\nu$ for $\nu \neq \nu_1, ..., \nu_j$. Then an integer $k = k(x) \in \{\nu_1, \nu_2, ..., \nu_j\}$ can be found such that $\tau_k := \min\{\tau_{\nu_1}, \tau_{\nu_2}, ..., \tau_{\nu_j}\}$. Define $u_x(\cdot) := u_{k(x)}(\cdot)$, $T(x) := \tau_{k(x)}$ and $\bar{q}(x) := q_{k(x)}$. Then, because of the local finiteness of covering, it follows by taking into account (2.9) that

$$\Phi(\pi(T(x), 0, x, u_x)) \leq \Phi(x) - \bar{q}(x), \tag{2.10a}$$

$$\Phi(\pi(t, 0, x, u_x)) \leq \bar{\beta}(\Phi(x)), \forall t \in [0, T(x)]; \forall x \in D_2 \setminus D_1 \tag{2.10b}$$

Furthermore, the mappings $T$ and $\bar{q}$ take strictly positive values, $T$ is LSC and satisfies (1.6) and $\bar{q}$ satisfies the following property:

$$\liminf_{x \to a} \bar{q}(x) > 0, \forall a \in D_2 \setminus D_1 \tag{2.11}$$

We are in a position to show (2.4). Define:

$$\varphi(t,s,x) := u_x(t-s), t \in [s, s+T(x)), x \in D_2 \setminus D_1 \tag{2.12}$$

where $u_x$ and $T$ are determined in (2.10). Let $x_0 \in D_2 \setminus N$ and consider the corresponding control $u_0 := u_{x_0}$ satisfying (2.10) with $x := x_0$. Define:

$$\begin{aligned} x_\nu &:= \pi(t_\nu, t_{\nu-1}, x_{\nu-1}, u_\nu), u_\nu(s) := \varphi(s, t_{\nu-1}, x_{\nu-1}) = u_{\nu-1}(s - t_{\nu-1}), s \in [t_{\nu-1}, t_\nu), \\ t_\nu - t_{\nu-1} &:= T(x_{\nu-1}), \nu = 1, 2, \ldots \text{ and } \Phi_\nu := \Phi(x_\nu), \nu = 0, 1, 2, \ldots \end{aligned} \tag{2.13}$$

provided that $x_\nu \in D_2 \setminus D_1$, $\nu = 1, 2, \ldots$, where $x_\nu|_{\nu=0}$ above, coincides with initial value $x_0$ and $t_0 := 0$. It follows from (2.10) and (2.13) that

$$\Phi_\nu \leq \Phi_{\nu-1} - \bar{q}(x_{\nu-1}) \tag{2.14a}$$
$$\Phi(\pi(s, t_{\nu-1}, x_{\nu-1}, u_\nu)) \leq \bar{\beta}(\Phi_{\nu-1}), s \in [t_{\nu-1}, t_\nu), \nu = 1, 2, \ldots \tag{2.14b}$$

therefore, by taking into account (2.14a), it follows:

$$\Phi_\nu \leq \Phi_0 - \sum_{i=0}^{\nu} \bar{q}(x_i), \nu = 1, 2, \ldots \tag{2.15}$$

From (2.2b,c),(2.13) and (2.15) we get

$$0 \leq \Phi_\nu \leq \Phi_0 < L_2, \nu = 1, 2, \ldots \tag{2.16}$$

Let $w = w_{x_0} : \mathbb{R}^+ \to \mathbb{R}^\ell$ be the concatenation of the controls $u_\nu : [t_{\nu-1}, t_\nu) \to \mathbb{R}^\ell$ above, namely:

$$w(t) := \varphi(t, t_{\nu-1}, x_{\nu-1}) = u_\nu(t), t \in [t_{\nu-1}, t_\nu) \tag{2.17}$$

Then by (2.14b),(2.16) and (2.17) it follows:

$$\Phi(\pi(s, t_{\nu-1}, x_{\nu-1}, w)) \leq \bar{\beta}(\Phi_{\nu-1}) \leq \bar{\beta}(\Phi_0), s \in [t_{\nu-1}, t_\nu), \nu = 1, 2, \ldots \tag{2.18}$$

provided that $x_\nu \in D_2 \setminus D_1$, $\nu = 1, 2, \ldots$. In order to establish that (2.4) holds for $x_0 \in D_2$ outside $N$, it suffices to show that

$$\exists \text{integer } \nu := \nu(x_0) : x_\nu \in \text{int } N \tag{2.19}$$

Indeed, notice first by virtue of (2.2b,c),(2.13) and (2.16) that

$$\{x_\nu, \nu = 1, 2, \ldots\} \subset D_2 \tag{2.20}$$

and suppose on the contrary that $\{x_\nu, \nu = 1, 2, \ldots\} \cap \text{int } N = \emptyset$. Then, there would exist a subsequence

$$\{x'_\nu\} \subset \{x_\nu\} \subset D_2 \setminus N \tag{2.21}$$

such that one of the following holds:

$$x'_\nu \to a \text{ for certain } a \in D_2 \cup \partial N \tag{2.22a}$$
$$x'_\nu \to a \text{ for certain } a \in \partial D_2 \tag{2.22b}$$

$$\lim_{\nu \to \infty} |x'_\nu| = \infty \qquad (2.22c)$$

To exclude (2.22a), we take into account (2.15),(2.16) and (2.20), which guarantee that $\sum_{i=0}^{\infty} \overline{q}(x_i) < \infty$. It then follows that $\overline{q}(x_\nu) \to 0$ and therefore $\overline{q}(x'_\nu) \to 0$ with $x'_\nu \to a \in D_2 \setminus D_1$, which contradicts (2.11). Suppose next that (2.22b) holds. Then by (2.2e) we may assume that without of generality it holds that $\Phi(x'_\nu) \to L_2$. On the other hand, by virtue of (2.16), $\limsup_{\nu \to \infty} \Phi(x'_\nu) \leq \Phi_0 < L_2$, a contradiction. Likewise, case (2.22c) is excluded, because of assumption (2.2e). Therefore, (2.4) is established and this completes the proof of statement (i).

In order to show statement (ii), we first use the result of statement (i), which, in our case (2.5) asserts that (2.4) holds for any neighborhood $N \subset \mathbb{R}^n$ of zero. The latter implies (1.5b). In order to show (1.5a), we invoke (2.5) and (2.2b,c,e), which guarantee existence of a function $a_1 \in K_\infty$ such that

$$a_1(|x|) \leq \Phi(x), \forall x \in D_2 \qquad (2.23)$$

It then follows, by virtue of (2.18),(2.20) and (2.23), that the corresponding trajectory $\pi(\cdot, 0, x_0; w)$ of (1.1a) satisfies:

$$|\pi(t, 0, x_0, w)| \leq a_1^{-1}(\overline{\beta}(\Phi(x_0))), \forall t \geq 0, x_0 \in D_2 \qquad (2.24)$$

which establishes (1.5a). The proof of last statement is an immediate consequence of statement (i) and (ii). The details are left to the reader.

The following remark, plays an important role in the proof of Proposition 3 in Section 3.

**Remark 2.** The result of statement (i) of Proposition 1 is also valid, if we replace (2.2e) plus (2.3) by the following hypothesis: There exist a continuous and positive definite function $\Phi: \mathbb{R}^n \to \mathbb{R}^+$ satisfying (2.2a,b,c), a function $\beta \in K$ and a map $q: \mathbb{R}^+ \to \mathbb{R}^+ \setminus \{0\}$ such that

$$\liminf_{s \to a} q(s) > 0, \forall a \in \mathbb{R}^+ \setminus \{0\}; \qquad (2.25)$$

$$\inf\{\Phi(x), x \in D_2 \setminus N\} > 0 \qquad (2.26)$$

for any closed neighborhood $N \subset D_2$ of $D_1$ and in such a way that, instead of (2.3), the following holds: for every constant $\sigma > 0$ and for every $x \in D_2 \setminus D_1$, a time $\tau = \tau(x) \in (0, \sigma]$ and a control $u_x(\cdot) : [0, \tau] \to \mathbb{R}^\ell$ can be found satisfying:

$$\Phi(\pi(\tau, 0, x, u_x)) < \Phi(x), \qquad (2.27a)$$

$$\Phi(\pi(t, 0, x, u_x)) \leq \beta(\Phi(x)) - q(\Phi(x)), \forall t \in [0, \tau] \qquad (2.27b)$$

The proof of the above claim is precisely the same with that given for the first statement of Proposition 1.

## 3. An Extension of Artstein-Sontag Theorem on Feedback Stabilization

The well known Artstein-Sontag theorem on stabilization for the affine in the control case (1.2) asserts that existence of a CLF, namely, of a $C^1$, positive definite and proper function $\Phi: \mathbb{R}^n \to \mathbb{R}^+$ satisfying the implication $(D\Phi g)(x) = 0, x \neq 0 \Rightarrow (D\Phi f)(x) < 0$, guarantees existence of a $C^\infty(\mathbb{R}^n \setminus \{0\})$ feedback law exhibiting global asymptotic stabilization for (1.2). We next apply the result of Section 2 to get a generalization of this result for the case (1.2). Our purpose is to show that, under weaker hypotheses, we may succeed SDF-WGAS for (1.2). Assume that $f$ and $g$ are $C^2$ and there exists a $C^2$, positive definite and proper function $\Phi: \mathbb{R}^n \to \mathbb{R}^+$ such that the following implication holds:

$$(g\Phi)(x) = 0, x \neq 0 \Rightarrow \begin{cases} \text{either } (f\Phi)(x) < 0, \\ \text{or } (f\Phi)(x) = 0; ([f,g]\Phi)(x) \neq 0 \end{cases} \quad (3.1)$$

where $[\cdot,\cdot]$ denotes the Lie bracket operator and we have used the standard notation $XY := DYX$ for any pair of $C^1$ mappings $X: \mathbb{R}^n \to \mathbb{R}^n, Y: \mathbb{R}^n \to \mathbb{R}^n$. We also recall the well known property that $[X,Y]\Phi = XY\Phi - YX\Phi$ for any pair of $C^2$ mappings $X,Y: \mathbb{R}^n \to \mathbb{R}^n$ and $C^2$ function $\Phi: \mathbb{R}^n \to \mathbb{R}$.

**Proposition 2.** Under previous hypotheses (3.1), system (1.2) is SDF-WGAS.

**Proof:** Let $\sigma > 0$ and state $x_0 \neq 0$. According to our assumption (3.1) we may distinguish three cases.
**Case 1:** $(g\Phi)(x_0) \neq 0$ and $(f\Phi)(x_0) < 0$.
Then, by considering the constant input

$$u_{x_0} := (-1 - (f\Phi)(x_0))/(g\Phi)(x_0) \quad (3.2)$$

there exists a time $\tau = \tau(x_0) \in (0, \sigma]$ such that the trajectory $\pi(\cdot) = \pi(\cdot, x_0, u) : [0, \tau] \to \mathbb{R}^n$ of system (1.2a) with $\pi(0, x_0, u) = x_0$ satisfies:

$$\Phi(\pi(t, x_0, u)) < \Phi(x_0), \forall t \in (0, \tau] \quad (3.3)$$

It turns out from (3.3) and Remark 1(i) that both (2.7a,b) hold with $\beta(s) = s$.
**Case 2:** $(g\Phi)(x_0) = 0$ and $(f\Phi)(x_0) < 0$.
Likewise for this case, we obtain (3.3), thus (2.7a,b) hold with $\beta(s) = s$ and zero input.
**Case 3:** $(g\Phi)(x_0) = 0$ and $(f\Phi)(x_0) = 0$.
Then, according to our assumption (3.1), it also holds:

$$([f,g]\Phi)(x_0) \neq 0 \quad (3.4)$$

Let $\tau \in (0, \tfrac{1}{2}\sigma]$ and a pair of constants $u_i, i = 1, 2$ yet to be determined. Let

$$X := f + u_2 g, \ Y := f + u_1 g \quad (3.5)$$

and consider the corresponding trajectories $X_t(r) := \pi(t, r, u_2)$, $Y_t(r) := \pi(t, r, u_1)$ of systems $\dot{x} = X(x)$ and $\dot{y} = Y(y)$, respectively, initiated at time $t = 0$ from $r \in \mathbb{R}^n$. Also, consider the trajectory of (1.2a):

$$R(t) := (X_t \circ Y_t)(x_0), \ t \in [0, \tfrac{1}{2}\sigma] \quad (3.6)$$

corresponding to the concatenation of $u_1$ and $u_2$ and of initial state $x_0$. Then we have:

$$\begin{aligned} \dot{R}(t) &= X(R(t)) + (DX_t Y) \circ Y_t(x_0) \\ &= X(R(t)) + ((DX_t Y) \circ X_{-t}) \circ R(t) \\ &= X(R(t)) + Y(R(t)) + \int_0^t (DX_s[Y,X]) \circ X_{-s}(R(t))ds \end{aligned} \quad (3.7)$$

and therefore:

$$\ddot{R}(t) = (X^2 + Y^2 + XY + YX) \circ R(t) + \left(DX(R(t)) + DY(R(t))\right)\int_0^t (DX_s[Y,X]) \circ X_{-s}(R(t))ds$$
$$+ \left(\int_0^t D(DX_s(R(t))[Y,X] \circ X_{-s}(R(t)))ds\right)\dot{R}(t) + (DX_t[Y,X]) \circ Y_t(x_0) \quad (3.8)$$

It follows from (3.7) and (3.8) that

$$\dot{R}(0) = X(x_0) + Y(x_0), \quad (3.9a)$$

$$\ddot{R}(0) = (X^2 + Y^2 + 2YX)(x_0) \tag{3.9b}$$

Since $(g\Phi)(x_0) = 0$ and $(f\Phi)(x_0) = 0$, it follows from (3.5) and (3.9a) that

$$D\Phi(x_0)\dot{R}(x_0) = D\Phi(x_0)(X(x_0) + Y(x_0)) = 0 \tag{3.10}$$

We next define:

$$m(t) := \Phi(R(t)) \tag{3.11}$$

Then, by invoking (3.5),(3.9), (3.10) and (3.11), we get:

$$\begin{aligned} m(t) &= m(0) + t\dot{m}(0) + \tfrac{t^2}{2}\ddot{m}(0) + 0(t^2) \\ &= \Phi(x_0) + t(X\Phi + Y\Phi)(x_0) + \tfrac{t^2}{2}\Big(\big((X+Y)^T D^2\Phi(X+Y)\big)(x_0) \\ &\quad + \big((X^2)\Phi + (Y^2)\Phi + 2(YX)\Phi\big)(x_0)\Big) + 0(t^2) \\ &= \Phi(x_0) + \tfrac{t^2}{2}\Big(\big((2f + (u_2 + u_1)g)^T D^2\Phi(2f + (u_2 + u_1)g)\big)(x_0) \\ &\quad + D\Phi(x_0)\big(4f^2 + (3u_2 + u_1)fg + (u_2 + 3u_1)gf + (u_2 + u_1)^2 g^2\big)(x_0)\Big) + 0(t^2) \end{aligned} \tag{3.12}$$

where $\lim_{t \to 0^+}(O(t)/t) = 0$. By defining $w := u_2 + u_1, u := u_1$, equality (3.12) is rewritten:

$$\begin{aligned} m(t) &= \Phi(x_0) + \tfrac{t^2}{2}\Big(\big((2f + wg)^T D^2\Phi(2f + wg)\big)(x_0) \\ &\quad + D\Phi(x_0)(4f^2 + (3w - 2u)fg + (w + 2u)gf + w^2 g^2)(x_0)\Big) + 0(t^2) \\ &= \Phi(x_0) + \tfrac{t^2}{2}\Big(\big((2f + wg)^T D^2\Phi(2f + wg)\big)(x_0) \\ &\quad + D\Phi(x_0)\big((4f^2)(x_0) + w(gf)(x_0) + 3w(fg)(x_0) \\ &\quad + 2u([g,f])(x_0) + w^2(g^2)(x_0)\big)\Big) + 0(t^2) \end{aligned} \tag{3.13}$$

By taking into account (3.4), it follows that every pair of constants $w$ and $c > 0$, a constant $u$ can be found such that

$$\begin{aligned} &\big((2f + wg)^T D^2\Phi(2f + wg)\big)(x_0) \\ &+ D\Phi(x_0)\big(4f^2(x_0) + w(gf)(x_0) + 3w(fg)(x_0)\big) \\ &+ 2u([g,f]\Phi)(x_0) + w^2(g^2)(x_0) \le -c \end{aligned} \tag{3.14}$$

and therefore from (3.13) and (3.14) we get:

$$\Phi(R(t)) < \Phi(x_0) - \tfrac{c}{2}t^2 + O(t^2) < \Phi(x_0), t > 0 \text{ near zero} \tag{3.15}$$

where $R(t) = (X_t \circ Y_t)(x_0)$; $X_t(r) = \pi(t, r, w - u)$, $Y_t(r) = \pi(t, r, u)$ with arbitrary constant input $w$ and constant control $u$ as defined in (3.14). As in the Cases 1 and 2, it follows from (3.15) that there exists a time $\tau = \tau(x_0) \in (0, \sigma]$ such that the trajectory $\pi(\cdot) = \pi(\cdot, x_0, u) : [0, \tau] \to \mathbb{R}^n$ of system (1.2a) satisfies (3.3), hence, we again conclude that (2.7a,b) hold with $\beta(s) = s$. It follows, by invoking Remark 1(i) and statement (iii) of Proposition 1, that system (1.2) is SDF-WGAS.

**Remark 3.** The result of Proposition 2 can directly be extended to multi-input affine in the control systems. Also, we note that, by considering higher derivatives of the map $t \to \Phi((X_t \circ Y_t)(x_0))$, we can get more general algebraic conditions in terms of $f$, $g$ and $\Phi$ for SDF-WGAS, provided that that dynamics of (1.2) and the function $\Phi$ are smooth enough.

**Example 1.** We illustrate the nature of Proposition 2 by considering the elementary planar case (1.2) with dynamics $f(x) := (a(x_1, x_2), a(x_1, x_2))^T$, $g(x) := (x_1, -x_2)^T$, $x := (x_1, x_2) \in \mathbb{R}^2$, where $a : \mathbb{R} \to \mathbb{R}^+$ is $C^1$ and satisfies $a(0,0) = 0$ and $x_1 a(x_1, x_1) < 0$ and $a(x_1, -x_1) \neq 0$ for every $x_1 \neq 0$. Then it can be easily verified that system above satisfies (3.1) with $V = \frac{1}{2}(x_1^2 + x_2^2)$, hence, it is SDF-WGAS.

## 4. SDF-WGAS for Composite Systems and the Small-Gain Property

Consider the composite system (1.3), where $f, g$ are Lipschitz continuous. For every $(x_0, y_0) \in \mathbb{R}^n \times \mathbb{R}^m$ we denote in the sequel by $(x(t), y(t)) = (x(t, x_0, y_0, u), y(t, x_0, y_0, u))$ the solution of (1.3) with initial $(x(0), y(0)) = (x_0, y_0) \in \mathbb{R}^n \times \mathbb{R}^m$ corresponding to certain control $u$. The main result of this section extends [34, Proposition 2].

For the subsystems $\Sigma_1, \Sigma_2$ in (1.3a) we make the following hypotheses:

**H1.** There exist continuous functions $V : \mathbb{R}^n \to \mathbb{R}^+$ and $a_1, a_2 \in K_\infty$ with

$$a_1(|x|) \leq V(x) \leq a_2(|x|), \quad \forall x \in \mathbb{R}^n \tag{4.1}$$

and functions $\beta_1, \gamma_1 \in K$ such that for any $\sigma > 0$ and for every $(x_0, y_0) \in \mathbb{R}^n \times \mathbb{R}^m$ with

$$|y_0| < \gamma_1(|x_0|) \tag{4.2}$$

there exist a time $\tau_1 = \tau_1(x_0, y_0) \in (0, \sigma]$ and a control $u_1 = u_{1(x_0, y_0)}(t) : [0, \tau_1] \to \mathbb{R}^\ell$ satisfying:

$$V(x(\tau_1, x_0, y_0, u_1)) < V(x_0) \tag{4.3a}$$
$$V(x(t, x_0, y_0, u_1)) \leq \beta_1(V(x_0)), \forall t \in [0, \tau_1] \tag{4.3b}$$
$$|y(t, x_0, y_0, u_1)| < \gamma_1(|x(t, x_0, y_0, u_1)|), \forall t \in [0, \tau_1] \tag{4.3c}$$

**H2.** There exist continuous functions $W : \mathbb{R}^m \to \mathbb{R}^+$ and $b_1, b_2 \in K_\infty$ with

$$b_1(|y|) \leq W(y) \leq b_2(|y|), \forall y \in \mathbb{R}^m \tag{4.4}$$

and functions $\beta_2, \Gamma_2 \in K$, such that for any $\sigma > 0$ and for every $(x_0, y_0) \in \mathbb{R}^n \times \mathbb{R}^m$ with

$$\Gamma_2(|x_0|) < |y_0| \tag{4.5}$$

there exist a time $\tau_2 = \tau_2(x_0, y_0) \in (0, \sigma]$ and a control $u_2 = u_{2(x_0, y_0)}(t) : [0, \tau_2] \to \mathbb{R}^\ell$ satisfying:

$$W(y(\tau_2, x_0, y_0, u_2)) < W(y_0) \tag{4.6a}$$
$$W(y(t, x_0, y_0, u_2)) \leq \beta_2(W(y_0)), \forall t \in [0, \tau_2] \tag{4.6b}$$
$$\Gamma_2(|x(t, x_0, y_0, u_2)|) < |y(t, x_0, y_0, u_2)|, \forall t \in [0, \tau_2] \tag{4.6c}$$

In addition to H1 and H2, the following hold:

**H3a.** *Small-Gain Property*: $\quad (b_1 \circ \gamma_1 \circ a_2^{-1})(s) > (b_2 \circ \Gamma_2 \circ a_1^{-1})(s), \forall s > 0$ (4.7)

**H3b.** For any $\sigma > 0$ and for every nonzero $(x_0, y_0) \in \mathbb{R}^n \times \mathbb{R}^m$ with

$$|y_0| < \gamma_1(|x_0|); \Gamma_2(|x_0|) < |y_0| \tag{4.8}$$

there exist a time $\tau = \tau(x_0, y_0) \in (0, \sigma]$ and a control $u = u_{(x_0, y_0)}(t) : [0, \tau] \to \mathbb{R}^\ell$ satisfying:

$$V(x(\tau, x_0, y_0, u)) < V(x_0) \tag{4.9a}$$
$$V(x(t, x_0, y_0, u)) \le \beta_1(V(x_0)), \forall t \in [0, \tau] \tag{4.9b}$$
$$W(y(\tau, x_0, y_0, u)) < W(y_0) \tag{4.9c}$$
$$W(y(t, x_0, y_0, u)) \le \beta_2(W(y_0)), \forall t \in [0, \tau] \tag{4.9d}$$

The main result in [34] guarantees that, under H1-H3, the composite system (1.3) is SDF-WGAS provided that $\gamma_1 \in K^\infty$. The following proposition establishes that, when $\gamma_1$ is bounded, system (1.3) is SDF-WGAS under H1-H3 plus some more appropriate hypotheses.

**Proposition 3.** For the subsystems system $\Sigma_1$, $\Sigma_2$ of (1.3) assume that H1, H2 and H3 are fulfilled. Then
(i) The composite system (1.3) is SDF-WGAS, provided that $\gamma_1 \in K^\infty$.
(ii) If $\gamma_1$ is bounded, then system (1.3) is SDF-WGAS, provided that, in addition to (4.7), it holds:

$$\lim_{s \to \infty}(b_1 \circ \gamma_1 \circ a_2^{-1})(s) > r := \lim_{s \to \infty}(b_2 \circ \Gamma_2 \circ a_1^{-1})(s) \tag{4.10}$$

and there exists a map $q : \mathbb{R}^+ \to \mathbb{R}^+ \setminus \{0\}$ satisfying (2.25) and such that, for any $\sigma > 0$ and for every $(x_0, y_0) \in \mathbb{R}^n \times \mathbb{R}^m$ with

$$W(y_0) \ge r \tag{4.11a}$$

there exist a time $\tau_2 = \tau_2(x_0, y_0) \in (0, \sigma]$ and a control $u_2 = u_{2(x_0, y_0)}(t) : [0, \tau_2] \to \mathbb{R}^\ell$ satisfying (4.6b) and simultaneously, instead of (4.6a), the following stronger property holds:

$$W(y(\tau_2, x_0, y_0, u_2)) < W(y_0) - q(W(y_0)) \tag{4.11b}$$

**Proof:** The proof of statement (i) has been given in [34]. Suppose next that $\gamma_1$ is bounded. Notice that (4.7) and boundedness of $\gamma_1$ implies that $\Gamma_2$ is also bounded and therefore $r < \infty$, as the latter defined in (4.10). To establish our claim we generalize the analysis adopted in [34]. First, by (4.7) and (4.10) we can find a pair of bounded functions $\ell_i \in K \cap C^1(\mathbb{R}^+ \setminus \{0\}), i = 1, 2$ with

$$(b_2 \circ \Gamma_2 \circ \alpha_1^{-1})(s) < \ell_1(s) < \ell_2(s) < (b_1 \circ \gamma_1 \circ \alpha_2^{-1})(s), \forall s > 0 \tag{4.12a}$$

and in such a way that, if we define:

$$R_i := \lim_{t \to \infty} \ell_i(t), i = 1, 2 \tag{4.12b}$$

then

$$R_1 > R_2 > r \tag{4.12c}$$

Obviously, since each $\ell_i$ is strictly increasing and bounded, we have:

$$R_i > \ell_i(s), \forall s \ge 0, i = 1, 2 \tag{4.13}$$

and the latter in conjunction with (4.1) implies:

$$\ell_i(V(x)) < R_i \text{ for all } x \in \mathbb{R}^n, i = 1, 2 \tag{4.14}$$

Next, by using (4.1), (4.4) and (4.12a), we can easily verify that for each $i = 1, 2$ it holds:

$$W(y_0) < \ell_i(V(x_0)), (x_0, y_0) \ne 0 \Rightarrow (4.2) \tag{4.15a}$$
$$W(y_0) > \ell_i(V(x_0)), (x_0, y_0) \ne 0 \Rightarrow (4.5) \tag{4.15b}$$

$$W(y_0) = \ell_i(V(x_0)), (x_0, y_0) \neq 0 \Rightarrow (4.8) \tag{4.15c}$$

Let $\bar{\beta}_1 \in K$ satisfying $(\ell_i \circ \beta_1)(s) \leq (\bar{\beta}_1 \circ \ell_i)(s)$, $\forall s > 0$, $i = 1, 2$. Then, by virtue of (4.3)-(4.6),(4.8),(4.9) and (4.15), the following properties hold for each $i = 1, 2$:

$$W(y_0) < \ell_1(V(x_0)), (x_0, y_0) \neq 0 \Rightarrow \begin{cases} \ell_i(V(x(\tau_1, x_0, y_0, u_1))) < \ell_i(V(x_0)) \\ \ell_i(V(x(t, x_0, y_0, u_1))) \leq \bar{\beta}_1(\ell_i(V(x_0))), \forall t \in [0, \tau_1] \\ W(y(t, x_0, y_0, u_1)) < \ell_i(V(x(t, x_0, y_0, u_1))), \forall t \in [0, \tau_1] \end{cases} \tag{4.16}$$

for certain $\tau_1 \in (0, \sigma]$ and control $u_1$;

$$\ell_i(V(x_0)) < W(y_0), (x_0, y_0) \neq 0 \Rightarrow \begin{cases} W(y(\tau_2, x_0, y_0, u_2)) < W(y_0) \\ W(y(t, x_0, y_0, u_2)) \leq \beta_2(W(y_0)), \forall t \in [0, \tau_2] \\ \ell_i(V(x(t, x_0, u_0, u_2))) < W(y(t, x_0, y_0, u_2)), \forall t \in [0, \tau_2] \end{cases} \tag{4.17}$$

for certain $\tau_2 \in (0, \sigma]$ with and control $u_2$;

$$\ell_i(V(x_0)) = W(y_0), (x_0, y_0) \neq 0 \Rightarrow \begin{cases} \ell_i(V(x(\tau, x_0, y_0, u))) < \ell_i(V(x_0)) \\ \ell_i(V(x(t, x_0, y_0, u))) \leq \beta_2(\ell_i(V(y_0))), \forall t \in [0, \tau] \\ W(y(\tau, x_0, u_0, u)) < W(y_0) \\ W(y(t, x_0, y_0, u)) \leq \beta_2(W(y_0)), \forall t \in [0, \tau] \end{cases} \tag{4.18}$$

for certain $\tau \in (0, \sigma]$ and control $u$.

Next, define $B(s) := 2\max\{\bar{\beta}_1(s), \beta_2(s)\}$ and consider the pair of continuous functions:

$$\Psi_i(x, y) := \begin{cases} W(y), & \text{for } W(y) \geq \ell_i(V(x)) \\ \ell_i(V(x)), & \text{for } \ell_i(V(x)) > W(y) \end{cases}, i = 1, 2 \tag{4.19}$$

By exploiting (4.16)-(4.19) for $i = 1$ and applying the same arguments with those in proof of Proposition 2 in [34], we can establish the following property:

**P1**. For each $\sigma > 0$ we can determine a map $T : (\mathbb{R}^n \times \mathbb{R}^m) \setminus \{0\} \to \mathbb{R}^+ \setminus \{0\}$ satisfying $T(x, y) \leq \sigma$ for all $(x, y) \in (\mathbb{R}^n \times \mathbb{R}^m) \setminus \{0\}$ and in such a way that for every $(x_0, y_0) \in (\mathbb{R}^n \times \mathbb{R}^m) \setminus \{0\}$, there exists a control $u = u_{(x_0, y_0)} : [0, T(x_0, y_0)] \to \mathbb{R}^\ell$ satisfying:

$$\Psi_1(x(T(x_0, y_0), x_0, y_0, u), y(T(x_0, y_0), x_0, y_0, u)) < \Psi_1(x_0, y_0) \tag{4.20a}$$
$$\Psi_1(x(t, x_0, y_0, u), y(t, x_0, y_0, u)) \leq B(\Psi_1(x_0, y_0)), \forall t \in [0, T(x_0, y_0)] \tag{4.20b}$$

For $i = 2$, we consider the function $\Psi_2(\cdot, \cdot)$ as defined by (4.19) and we exploit our additional assumption that (4.11a) implies both (4.6b) and (4.11b). Notice that, when $W(y) > R_2$, it follows from (4.12c) and (4.14) that $W(y) > r > \ell_2(V(x))$ thus, by (4.19) we have $\Psi_2(x, y) = W(y)$. The latter in conjunction with second inequality in (4.17) for $i = 2$, (which is a consequence of (4.6b)), together with implication (4.11a) $\Rightarrow$ (4.11b) and definition of $B(\cdot)$, imply the following property:

**P2**. For every $\sigma > 0$ there exists a map $T : (\mathbb{R}^n \times \mathbb{R}^m) \setminus \{0\} \to \mathbb{R}^+ \setminus \{0\}$ satisfying $T(x, y) \leq \sigma$ for all $(x, y) \in (\mathbb{R}^n \times \mathbb{R}^m) \setminus \{0\}$ such that for every $(x_0, y_0) \in (\mathbb{R}^n \times \mathbb{R}^m) \setminus \{0\}$ with $W(y_0) \geq R_2$, a control $u = u_{(x_0, y_0)} : [0, T(x_0, y_0)] \to \mathbb{R}^\ell$ can be found satisfying:

$$\Psi_2(x(T(x_0, y_0), x_0, y_0, u), y(T(x_0, y_0), x_0, y_0, u)) < \Psi_2(x_0, y_0) - q(\Psi_2(x_0, y_0)); \tag{4.21a}$$
$$\Psi_2(x(t, x_0, y_0, u), y(t, x_0, y_0, u)) \leq B(\Psi_2(x_0, y_0)), \forall t \in [0, T(x_0, y_0)]; \tag{4.21b}$$

We are in a position, by exploiting Properties P1 and P2 and the results of Proposition 1 and Remark 2, to establish that (1.3) is SDF-WGAS. We distinguish two cases:

**Case 1:** $(x_0, y_0) \in \mathbb{R}^n \times \mathbb{R}^m : W(y_0) < R_1$, as the latter is defined by (4.12b).
We define:

$$\Phi(x, y) := \Psi_1(x, y), \ \beta := B, \ L_1 := 0, L_2 := R_1; \qquad (4.22a)$$

$$D_1 := \{(0,0) \in \mathbb{R}^n \times \mathbb{R}^m\}; \qquad (4.22b)$$

$$D_2 := \{(x, y) \in \mathbb{R}^n \times \mathbb{R}^m : R_1 > W(y) \geq 0\} \qquad (4.22c)$$

Then, we can verify that all conditions of second statement of Proposition 1 are fulfilled with $\Phi, \beta, L_1, L_2, D_1$ and $D_2$ as above. Indeed, note that, for $(x, y) \in D_2 \setminus \{0\}$, we have, due to (4.19) and (4.22a,c), that $\Psi_1(x, y) = W(y) \leq R_1 (= L_2)$ for the case $\ell_1(V(x)) \leq W(y)$ and, due to (4.14),(4.19) and (4.22a,c), it holds that $\Psi_1(x, y) = \ell_1(V(x)) \leq R_1 (= L_2)$ when $\ell_1(V(x)) \geq W(y)$, hence, (2.2b) is satisfied. Also, for $(x, y) \notin D_2$ we have, due to (4.14), that $W(y) \geq R_1 > \ell_1(V(x))$, hence, by (4.19) it follows that $\Psi_1(x, y) = W(y) \geq R_1 (= L_2)$, therefore (2.2c) holds as well. We next show (2.2e). Obviously, by (4.22c) we have $W(a_2) = R_1$ for every $a := (a_1, a_2) \in \partial D_2$, which, in conjunction with (4.14) and (4.19), implies that $\Psi_1(a) = W(a_2) = R_1$. Consider next a sequence $\{(x_\nu, y_\nu) \in D_2\}$ with $\lim_{\nu \to \infty} |(x_\nu, y_\nu)| = \infty$. We show that $\limsup_{x_\nu \to \infty} \Psi_1(x_\nu, y_\nu) = R_1$. Because of (4.4) and (4.22c), the sequence $\{y_\nu\}$ is bounded, hence, $\lim_{\nu \to \infty} |x_\nu| = \infty$. Then, by (4.1) and (4.12b) we have:

$$\lim_{x_\nu \to \infty} \ell_i(V(x_\nu)) = R_1 \qquad (4.23)$$

Since $R_1 > W(y_\nu), \nu = 1, 2, \ldots$, we may distinguish two cases. The first is $R_1 > W(y'_\nu) \geq \ell_1(V(x'_\nu))$ for certain subsequence $\{(x'_\nu, y'_\nu)\} \subset \{(x_\nu, y_\nu)\}$ and therefore by (4.19), $R_1 > \Psi_1(x'_\nu, y'_\nu) \geq \ell_1(V(x'_\nu))$. The latter in conjunction with (4.23) implies:

$$\lim_{\nu \to \infty} \Psi_1(x'_\nu, y'_\nu) = R_1 \qquad (4.24)$$

The other case is $R_1 > \ell_1(V(x'_\nu)) > W(y'_\nu)$ for certain subsequence $\{(x'_\nu, y'_\nu)\} \subset \{(x_\nu, y_\nu)\}$. Then by (4.19) it follows that $\Psi_1(x'_\nu, y'_\nu) = \ell_1(V(x'_\nu))$, hence, by (4.23) we again obtain (4.24). The above discussion asserts that (2.2e) is fulfilled. Finally, (2.3a,b) is a consequence of Property P1 and definitions (4.22). We can therefore apply the result of Proposition 1(ii) to establish existence of a sampled-data controller $u = u_{(x_0, y_0)}(\cdot)$, $(x_0, y_0) \in D_2$, such that the corresponding trajectory $(x(\cdot, x_0, y_0, u_{(x_0, y_0)}), y(\cdot, x_0, y_0, u_{(x_0, y_0)}))$ of (1.3a) is tending to zero $(0,0) \in \mathbb{R}^n \times \mathbb{R}^m$ as $t \to \infty$. Also, according to the second statement of Proposition 1, we have:

$$\forall \varepsilon > 0 \Rightarrow \exists \delta = \delta(\varepsilon) > 0:$$
$$|(x_0, y_0)| \leq \delta \Rightarrow |(x(t, x_0, y_0, u_{(x_0, y_0)}), y(t, x_0, y_0, u_{(x_0, y_0)}))| \leq \varepsilon, \ \forall t \geq 0 \qquad (4.25)$$

**Case 2:** $(x_0, y_0) \in (\mathbb{R}^n \times \mathbb{R}^m) : W(y_0) \geq R_1$.
Let $\bar{R}$ be a constant such that $R_1 > \bar{R} > R_2$ and define:

$$\Phi(x, y) := \Psi_2(x, y), \ \beta := B, \ L_1 := R_2, L_2 := \infty; \qquad (4.26a)$$

$$D_1 := \{(x, y) \in \mathbb{R}^n \times \mathbb{R}^m : W(y) \leq R_2\}; \qquad (4.26b)$$

$$N := \{(x, y) \in \mathbb{R}^n \times \mathbb{R}^m : W(y) \leq \bar{R}\}; \qquad (4.26c)$$

$$D_2 := \mathbb{R}^n \times \mathbb{R}^m; \qquad (4.26d)$$

According to definitions (4.26) and by taking into account Property P2, we can easily verify that all conditions of Remark 2 are satisfied with $\Phi$, $\beta$, $L_1, L_2$ $D_1$ and $D_2$ as above. (For completeness, we note that (2.26) and (2.27) are consequence of (4.4),(4.26),(4.19) and (4.21)). Therefore, according to the Remark 2, there exists a sampled-data controller $u = u_{(x_0, y_0)}(\cdot), (x_0, y_0) \in (\mathbb{R}^n \times \mathbb{R}^m) \setminus D_1$ such that for $(x_0, y_0) \in (\mathbb{R}^n \times \mathbb{R}^m) \setminus N$ the corresponding trajectory $(x(\cdot, x_0, y_0, u_{(x_0, y_0)}), y(\cdot, x_0, y_0, u_{(x_0, y_0)}))$ of (1.3a) enters int $N$ after some finite time, say $\tau_1$; particularly, if we denote:

$$(x_1, y_1) := (x(\tau_1, x_0, y_0, u_{(x_0, y_0)}), y(\tau_1, x_0, y_0, u_{(x_0, y_0)}))  \qquad (4.27)$$

then we have:

$$\Psi_2(x_1, y_1) = W(y_1) < \bar{R} \qquad (4.28)$$

therefore $W(y_1) < R_1$. We then apply the arguments used for the Case 1 with $(x_0, y_0) := (x_1, y_1)$ as above and find a sampled-data controller $\bar{u} = \bar{u}_{(x_1, y_1)}(\cdot)$ such that the trajectory $(x(\cdot, x_1, y_1, \bar{u}_{(x_1, y_1)}), y(\cdot, x_1, y_1, \bar{u}_{(x_1, y_1)}))$ of (1.3a) is tending to zero $(0,0) \in \mathbb{R}^n \times \mathbb{R}^m$ as $t \to \infty$. It turns out, by considering the concatenation $w = w_{(x_0, y_0)}(\cdot)$ of the sampled-data controls $u = u_{(x_0, y_0)}(\cdot)$, $t \in [0, \tau_1)$ and $\bar{u} = \bar{u}_{(x_1, y_1)}(\cdot), t \in [\tau_1, \infty)$, that the corresponding trajectory $(x(\cdot, x_1, y_1, w_{(x_0, y_0)}), y(\cdot, x_1, y_1, w_{(x_0, y_0)}))$ of (1.3a) is tending to zero $(0,0) \in \mathbb{R}^n \times \mathbb{R}^m$ as $t \to \infty$. The latter in conjunction with (4.25) asserts that the composite system (1.3) is SDF-WGAS.

**Remark 4.** (i) It can be easily verified that implication (4.11a) $\Rightarrow$ ((4.6b) plus (4.11b)) is fulfilled, if $W$ is $C^1$ and there exists a function $d \in K$, such that for every $(x, y) \in \mathbb{R}^n \times \mathbb{R}^m$ with $W(y) \geq r$ there exists a vector $u \in \mathbb{R}^m$ such that $DW(y)G(x, y, u) \leq -d(W(y))$.

(ii) As is pointed out in [34], the sufficient conditions of Proposition 3 are simplified when, either $f$, or $g$ is independent of $u$. Particularly, for the case of systems $\Sigma_1 : \dot{x} = f(x, y)$; $\Sigma_2 : \dot{y} = g(x, y, u)$ assume that, instead of H1 for $\Sigma_1$, there exist continuous functions $V : \mathbb{R}^n \to \mathbb{R}^+$ and $a_1, a_2 \in K_\infty$, $\gamma_1, \beta_1 \in K$ such that (4.1) holds and in such a way that for every $(x_0, y_0) \in \mathbb{R}^n \times \mathbb{R}^m, y_0 \neq 0$ for which (4.2) holds and for any $\sigma > 0$ there exists a time $\tau_1 \in (0, \sigma]$ satisfying

$$V(x(\tau_1, x_0, y_0)) < V(x_0) \qquad (4.29a)$$
$$V(x(t, x_0, y_0)) \leq \beta_1(V(x_0)), \forall t \in [0, \tau_1] \qquad (4.29b)$$

where $x(\cdot, x_0, y_0)$ denotes the solution of $\Sigma_1$. Also, assume that $\Sigma_2$ satisfies H2 and, either $\gamma_1$ is unbounded, or property (4.10) and implication (4.11a) $\Rightarrow$ (4.11b) plus (4.6b) are fulfilled for the case where $\gamma_1$ is bounded. Then, under H3a the system above is SDF-WGAS by means of a static sampled-data feedback stabilizer; (the additional condition H3b is not required). Implication (4.2) $\Rightarrow$ (4.29a,b) is fulfilled with $\beta_1(s) = s$, if we assume that $\Sigma_1$ satisfies the Lyapunov characterization of the version of the ISS property, with input $y$ and state $x$, which is originally introduced in [31].

As a consequence of Proposition 3(i) we provide the following result concerning the case (1.4):

**Proposition 4.** For system (1.4) we assume that $A_i, B_i$ are $C^\infty$ and there exist functions $\gamma_1, \Gamma_2 \in K^\infty$ with

$$\gamma_1(s) > \Gamma_2(s), \forall s > 0 \qquad (4.30)$$

such that

- $\dim Lie\{A_i(\bullet, y), i = 1, ..., \ell\}|_x = n$, $\forall (x, y) \in \mathbb{R}^n \times \mathbb{R}^m$, $x \neq 0$ with $|y| < \gamma_1(|x|)$ (4.31)
- $\dim Lie\{B_i(x, \bullet), i = 1, ..., \ell\}|_y = m$, $(x, y) \in \mathbb{R}^n \times \mathbb{R}^m$, $y \neq 0$ with $|y| > \Gamma_2(|x|)$ (4.32)

- $$\dim Lie\{F_i, i=1,...,\ell\}|_\xi = n+m,$$
$$\forall \xi = (x,y) \in \mathbb{R}^n \times \mathbb{R}^m \text{ with } |y| < \gamma_1(|x|) \text{ and } |y| > \Gamma_2(|x|) \tag{4.33}$$

Then (1.4) is SDF-WGAS.

**Remark 4.** If the assumptions of Proposition 4 are substituted by the stronger accessibility rank condition for the whole system (1.4a,b), namely, $\dim Lie\{F_i, i=1,...,\ell\}|_\xi = n+m$ for every nonzero $\xi = (x,y) \in \mathbb{R}^n \times \mathbb{R}^m$, then the result of Proposition 4 is strengthened; particularly, according to the well known result in [9], this condition guarantees that the system (1.4) is *globally stabilizable by means of a smooth time-varying feedback*. A variant of this condition has been used in [33, Corollary 1] to establish SDF-WGAS for the nonholonomic case (1.4).

**Proof of Proposition 4:** Consider the case $(x_0, y_0) \in \mathbb{R}^n \times \mathbb{R}^m$ with

$$|y_0| < \gamma_1(|x_0|) \tag{4.34}$$

and for the given $y_0$ above consider the auxiliary system:

$$\dot{X} = f(X, y_0, u) := \sum_{i=1}^{\ell} u_i A_i(X, y_0, u) \tag{4.35a}$$

$$\dot{Y} = g(X, y_0, u) := \sum_{i=1}^{\ell} u_i B_i(X, y_0, u) \tag{4.35b}$$

and let us denote by $(X(\cdot), Y(\cdot)) = (X(\cdot, x_0, y_0, u), Y(\cdot, x_0, y_0, u))$ its trajectory initiated from $(x_0, y_0) \in \mathbb{R}^n \times \mathbb{R}^m$ at time $t=0$ and corresponding to some control $u$. Also define $V(x) := |x|^2$. Then using (4.31) and applying the result of Corollary 1 in [33] for the subsystem (4.35a) it follows there exist $\beta_1 \in K$ such that for any $\sigma > 0$ there exists a time $\tau_1 \le \sigma$ and a control $u_1 = u_{1(x_0, y_0)}(t) : [0, \tau_1] \to \mathbb{R}^\ell$ such that for $(x_0, y_0) \in \mathbb{R}^n \times \mathbb{R}^m$ satisfying (4.34), the corresponding trajectory

$$(X(\cdot), Y(\cdot)) = (X(\cdot, x_0, y_0, u_1), Y(\cdot, x_0, y_0, u_1)) : [0, \tau_1] \to \mathbb{R}^\ell$$

of (4.35) satisfies:

$$V(X(\tau_1, x_0, y_0, u_1)) < V(x_0) \tag{4.36a}$$
$$V(X(t, x_0, y_0, u_1)) \le \beta_1(V(x_0)), \forall t \in [0, \tau_1] \tag{4.36b}$$
$$|Y(t, x_0, y_0, u_1)| < \gamma_1(|X(t, x_0, y_0, u_1)|), \forall t \in [0, \tau_1] \tag{4.36c}$$

It turns out, by taking into account (4.36) and using elementary continuity arguments, that there exist a time $\tau_1' \le \tau_1 (\le \sigma)$ such that

$$V(x(\tau_1', x_0, y_0, u_1)) < V(x_0) \tag{4.37a}$$
$$V(x(t, x_0, y_0, u_1)) \le 2\beta_1(V(x_0)), \forall t \in [0, \tau_1'] \tag{4.37b}$$
$$|y(t, x_0, y_0, u_1)| < \gamma_1(|x(t, x_0, y_0, u_1)|), \forall t \in [0, \tau_1'] \tag{4.37c}$$

where $(x(\cdot), y(\cdot)) = (x(\cdot, x_0, y_0, u), y(\cdot, x_0, y_0, u))$ denotes the solution of the original system (1.4a,b) with initial $(x(0), y(0)) = (x_0, y_0)$. We conclude that Hypotheses H1 of Proposition 3 is satisfied with $V(x) := |x|^2$. Likewise, for the case $(x_0, y_0) \in \mathbb{R}^n \times \mathbb{R}^m$ with $|y_0| > \Gamma_2(|x_0|)$ we consider for this $x_0$ the auxiliary system:

$$\dot{X} = f(x_0, Y, u) := \sum_{i=1}^{\ell} u_i A_i(x_0, Y, u) \tag{4.38a}$$

$$\dot{Y} = g(x_0, Y, u) := \sum_{i=1}^{\ell} u_i B_i(x_0, Y, u) \tag{4.38b}$$

Then by invoking (4.32), defining $W(y) := |y|^2$ and by again applying Corollary 1 in [33] for the subsystem (4.38b) we can verify as above that Hypothesis H2 of Proposition 3 is satisfied as well. Likewise, by exploiting (4.33), we can verify hat H3a and H3b hold for (1.4) with $V$ and $W$ as above, $a_1 = a_2 \equiv s^2$, $b_1 = b_2 \equiv s^2$ and $\gamma_1, \Gamma_2 \in K$ as given in our statement. We conclude, according to the result of Proposition 3, that system (1.4) is SDF-WGAS.

**Example 2.** We illustrate the nature of Proposition 4 by considering the planar case:

$$\begin{pmatrix} \dot{x} \\ \dot{y} \end{pmatrix} = u_1 \begin{pmatrix} A_1(x,y) \\ B_1(x,y) \end{pmatrix} + u_2 \begin{pmatrix} A_2(x,y) \\ B_2(x,y) \end{pmatrix}, (x,y) \in \mathbb{R}^2 \quad (4.39)$$

where $A_i, B_i$ are $C^1$ and satisfy:

$$A_1(0,\cdot) = A_2(0,\cdot) = 0, B_1(\cdot,0) = A_2(\cdot,0) = 0 \quad (4.40)$$

$$(A_1)_x A_2 \big|_{(x,y)} \neq (A_2)_x A_1 \big|_{(x,y)}, \forall (x,y) \in \mathbb{R}^2, x \neq 0 \quad (4.41a)$$

$$(B_1)_y B_2 \big|_{(x,y)} \neq (B_2)_y B_1 \big|_{(x,y)}, \forall (x,y) \in \mathbb{R}^2, y \neq 0 \quad (4.41b)$$

$$A_1 B_2 \big|_{(x,y)} \neq A_2 B_1 \big|_{(x,y)}, \forall (x,y) \in \mathbb{R}^2, x \neq 0, y \neq 0 \quad (4.42)$$

By taking arbitrary $\gamma_1, \Gamma_2 \in K^\infty$ satisfying (4.30), we can easily verify that all conditions of Proposition 4 are fulfilled, therefore (4.39) is SDF-WGAS. For completeness we note that (4.41a), (4.41b) imply (4.31) and (4.32), respectively and (4.42) implies (4.33).